\documentclass[12pt]{amsart}

\usepackage{amssymb}

\usepackage{enumitem}
\usepackage{mathbbol}
\usepackage{graphicx}

\makeatletter
\@namedef{subjclassname@2020}{%
  \textup{2020} Mathematics Subject Classification}
\makeatother

\usepackage[T1]{fontenc}

\newtheorem{theorem}{Theorem}[section]

\newtheorem{lemma}[theorem]{Lemma}

\theoremstyle{definition}
\newtheorem{definition}[theorem]{Definition}
\newtheorem{remark}[theorem]{Remark}

\numberwithin{equation}{section}

\begin{document}


\baselineskip=17pt


\title[A Discussion of Arnold's Limit Problem and its Geometric Argument]{A Discussion of Arnold's Limit Problem and its Geometric Argument}

\author[K. Honn]{Keising Honn}
\address{School of Mathematics and Statistics, Lanzhou University, No. 222 South Tianshui Road, Lanzhou 730000, Gansu Province, P.R.China}
\email{hanqch21@lzu.edu.cn}

\date{}

\begin{abstract}
	Upon re-examining Arnold's established lemma for explaining his famous limit problem, we have determined that while the lemma itself is correct, there is a defect in the original geometric proof. In this paper, we prove the correctness of the lemma using methods of power series, and construct a counterexample to illustrate the defect in Arnold's geometric proof.
\end{abstract}



\maketitle

\section{Introduction}~\

\begin{figure}[h!]
	\centering
	\includegraphics{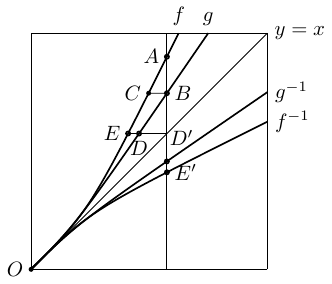}
	\caption{$f$ and $f^{-1}$, $g$ and $g^{-1}$ are symmetric about $y=x$ respectively}
	\label{fig: 1}
\end{figure}

Arnold's limit problem is
$$\lim_{x \rightarrow 0} \frac{\tan \sin x - \sin \tan x}{\arctan \arcsin x - \arcsin \arctan x} = ?$$
In general form, this problem can be formulated as
$$\lim_{x \rightarrow 0} \frac{f(x) - g(x)}{g^{-1}(x) - f^{-1}(x)} = ? \eqno(1)$$
with
$$f,\,g \text{ analytic}, \quad f(0) = g(0) = 0, \quad f'(0) = g'(0) = 1.\eqno(2)$$
In \cite{Arn90}, Arnold presented his famous geometric solution. The core of his method is
\begin{lemma}\label{lemma: 1}
	If the graphs of analytic functions $f$ and $g$ do not coincide and they are both tangent to the line $y=x$ at the origin (Fig.\ref{fig: 1}), then $|AB|/|BC|$ and $|BC|/|ED|$ converges to 1 as $A$ is sufficiently close to the origin.
\end{lemma}

Arnold's own proof of Lemma \ref{lemma: 1} is Fig.\ref{fig: 1}. Intuitively, it seems that when $A$ close to $O$, $AC$ as a cut line transitions to a tangent line, in turn, $ACB$ tends to an isosceles right triangle and $BCED$ tends to a parallelogram. Hence limit (1) is equivalent to
$$\lim \frac{|AB|}{|D'E'|} = \lim \frac{|AB|}{|ED|} = \lim \left(\frac{|AB|}{|BC|} \cdot \frac{|BC|}{|ED|}\right) = 1.$$
There is no analytic formulation of this Lemma (and its proof above) has appeared in the current literature. Therefore, we first solve problem (1) -- (2) by means of power series, which shows the correctness of Lemma \ref{lemma: 1} itself. However, when trying to describe ``$BCED$ tends to a parallelogram'' in analysis, we found a defect in the above geometric proof. Furthermore, a counterexample is constructed to confirm the nature of this defect.

\section{An Analysis Argument of Problem (1) -- (2)}~\

Since $f$ and $g$ are analytic, they can be represented as power series respectively. We denote by
$$f(x) = \sum_{k \geq 1} a_k x^k \quad \text{and} \quad f^{-1}(x) = \sum_{l \geq 1} b_l x^l.\eqno(3)$$
The summations starts with $k,\,l=1$, since by the restriction $f(0) = 0$, the coefficient of $x^0$ must be 0; moreover, due to the existence of inverse there must be $a_1,\,b_1 \neq 0$.\par
The connection between two power series in (3) is inscribed by the following Theorem.
\begin{theorem}
	For coefficients in (3) we have
	$$
		b_n = \left\{
			\begin{aligned}
			&1/a_1,\quad &n = 1,\\
			&-a_n/a_1^{n+1}+R_n,\quad &n > 1,
			\end{aligned}
			\right.$$
	where $R_n$ are some terms determined only by $a_k$ for $k < n$.
\end{theorem}
\begin{proof}
	Note that
$$f \circ f^{-1}(x) = \sum_{k\geq 1} a_k\left(\sum_{l\geq 1} b_l x^l\right)^k = x := \sum_{n \geq 1} c_n x^n,\eqno(4)$$
here the rightmost end is clearly a power series of the form, with all coefficients 0 except $c_1 = 1$.
	For each $n > 1$, $c_n$ does not absorb $b_l$ that makes $l > n$, otherwise the power of $x$ in such $b_l$-term will greater than $n$. Similarly $c_n$ does not contain those $a_k$ that makes $k > n$. However, $\sum_{l \geq 1} b_l x^l$ already contains $b_n x^n$, which means that $k$ here can only be taken as 1, so the $b_n$-term involved can only be $a_1b_nx^n$. Similarly, the $a_n$-term in $a_n \left(\sum_{l \geq 1} b_l x^l\right)^n$ can only be $a_n b_1^n x^n$. Therefore,
$$c_n = 0 = a_1 b_n + a_n b_1^n + K,$$
where $K$ is some terms determined only by $a_k$ and $b_l$ for $l,\,k < n$; and for $c_1 = 1$ we obviously have $c_1 = a_1 b_1$. Thus we have shown that, $b_n$ can only be determined by $a_k$ for $k \leq n$ and $b_l$ for $l < n$.\par
Furthermore, by $b_1 = 1/a_1$, $b_1$ is completely determined by $a_1$. $b_2$ is determined by $a_1,\,a_2,\,b_1$, which means by $a_1,\,a_2$. Continuing this process, eventually $b_n$ is completely determined by $a_k$ only for $k \leq n$, which finishes our proof.
\end{proof}
\begin{remark}
We have in fact verified an even stronger conclusion, which is that, $R_n$ is the same for all $f$ satisfying condition (2).
\end{remark}
Fix
$$g(x) = \sum_{k \geq 1} A_k x^k \quad \text{and} \quad g^{-1}(x) = \sum_{l \geq 1} B_l x^l,$$
then since $f'(0) = g'(0) = 1$, we have $a_1 = A_1 = 1$. Hence
$$b_n = a_n + R_n,\quad B_n = -A_n + R_n.$$
It is clear that $b_n = B_n$ holds only if $a_n = A_n$. Suppose that the first distinct coefficients in the power series of $f$ and $g$ occur in the place of $N$-th term. In this situation, $b_n - B_N = -(a_N - A_N)$, thus
$$\lim_{x \rightarrow 0} \frac{f(x) - g(x)}{g^{-1}(x) - f^{-1}(x)} = \lim_{x \rightarrow 0} \frac{(a_N - A_N)x^N + O(x^{N+1})}{(a_N - A_N)x^N + O(x^{N+1})}  = 1.$$

\section{Motivation to Detect a Defect}~\

Recall Lemma \ref{lemma: 1} and Fig.\ref{fig: 1}, we will notice that the Lemma itself requires the assumption of $f,\,g$ analytic, but this condition is not used in Arnold's geometric proof, because we cannot represent such a strong condition of analyticity just in a graph; at least the graph of an analytic function does not look different from the graph of a ``smooth'' function. So this geometric proof may be tempting to think that Lemma \ref{lemma: 1} may still holds, after weakening the original assumption to ``$f,\,g$ are `smooth' functions (e.g. $C^\infty$-functions)''. Under the latter weaker condition we try to provide an analysis formulation of limit processes appearing in the geometric proof of Section 1, and the defect appears. To this end, without loss of generality, we suppose the coordinates (see Fig.\ref{fig: 1})
$$\begin{aligned}
	A(x,f(x)),\quad B(x,\,g(x)),\quad C(f^{-1}\circ g(x),g(x)),\\
	D(g^{-1}(x), x),\quad D'(x,x),\quad E(f^{-1}(x),x).
\end{aligned}$$
\subsection{Analysis argument of $|AB|/|BC| \rightarrow 1$. }We have
$$|AB| = f(x) - g(x),\quad |BC| = x - f^{-1} \circ g(x),$$
which yields
$$\frac{|AB|}{|BC|} = \frac{f(x) - f(f^{-1}\circ g(x))}{x - f^{-1}\circ g(x)}.\eqno(5)$$
Applying the Midian Theorem to (5), we can find a point $x \leq \xi \leq f^{-1}\circ g(x)$ such that $|AB| / |BC| = f'(\xi)$. But $f^{-1}(x) \circ g(x) \rightarrow 0$ as $x \rightarrow 0$, so we have $f'(\xi) \rightarrow f'(0) = 1$ by continuity.
\subsection{Analysis argument of $|ED|/|BC| \rightarrow 1$. }
\begin{figure}[h!]
        \centering
        \begin{minipage}{0.45\linewidth}
                \centering
		\includegraphics{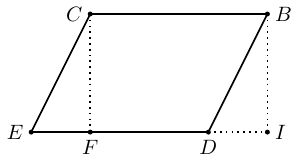}
                \caption{Convergence}
                \label{fig:2}
        \end{minipage}
\begin{minipage}{0.45\linewidth}
                \centering
		\includegraphics{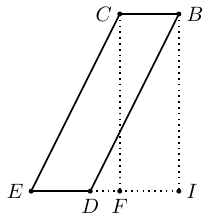}
                \caption{Divergence}
                \label{fig:3}
        \end{minipage}
\end{figure}

The defect in Arnold's geometric proof will be cited here. In this case we should note that, $BCED$ is never a parallelogram, but tends to be as $A$ approaches $O$. So as shown in Fig.\ref{fig:2}, 
$$|ED| = |EF| + |FD| = |DD'|\frac{|EF|}{|DD'|} + |FD| = |FD'| + |DD'| \left(\frac{|EF|}{|DD'|} - 1\right).$$
Note that $|BC| \sim |FD'|$ as $A \rightarrow O$, we have
$$\frac{|ED|}{|BC|} \sim 1 + \frac{|DD'|}{|FD'|}\left(\frac{|EF|}{|DD'|} - 1\right).\eqno(6)$$
In the case of Fig.\ref{fig:2}, we have $|EF|/|DD'| - 1 \rightarrow 0$ and $|DD'| / |FD'|$ bounded (since $|DD'| < |FD'|$), so there is indeed $|ED|/|BC| \rightarrow 1$. But if we consider case of Fig.\ref{fig:3}, i.e., the quasi-parallelogram is very narrow and $|BC|,\,|DE|$ are small enough, it is easy to verify that (6) still holds, but $|DD'| > |FD'|$, which is likely to lead to $|DD'|/|FD'|$ diverges. This situation will always occur if $f$ and $g$ are moving closer together than they are to $y = x$. This suggests to us that there exists such $C^\infty$-functions, which make Lemma \ref{lemma: 1} incorrect, but Arnold's geometric proof still holds.

\section{Construction of the Counterexample}~\

We start with
\begin{definition}
	We define $\theta(x): \mathbb R \rightarrow \mathbb R$ by
$$\theta(x) = \left\{
	\begin{aligned}
		e^{-1/|x|},&\quad x \neq 0,\\
		0,&\quad x = 0.
	\end{aligned}
	\right.$$
\end{definition}
	\begin{theorem}
		$\theta(x) \in C^\infty$ but not analytic.
	\end{theorem}
	\begin{proof}
		We have
		$$\theta(0) = \theta'(0) = \theta ''(0) = \cdots = \theta^{(n)}(0) = 0$$
		but for any $x,\,\xi > 0$ we have $\theta(x) \neq 0$ and
		$$\lim_{n \rightarrow \infty} \sum_{k=1}^n \frac{\theta^{(k)}(\xi)}{k!} (x - \xi)^k = 0.$$
		Note that $\theta$ is even, we finishes our proof.
	\end{proof}
	The whole problem is symmetric for $f,\,g$ and their inverses. Define $p = f^{-1}$ and $q = g^{-1}$, then
	$$\frac{|BC|}{|ED|} = \frac{x - p\circ g(x)}{q(x) - p(x)}.$$
	In order for $p,\,q$ to also satisfy the requirements stated at the end of Section 3, we may set
	$$p(x) = q(x) + \theta(x), \quad q(x) = x + x^2,\eqno(7)$$
	we choose the $x^2$ term in (7) so that it serves the purpose of controlling graphs of $p,\,q$ away from $y=x$. Denote $t = g(x)$, hence $x = q(t)$, and
	$$\frac{|BC|}{|ED|} = \frac{q(t) - p(t)}{\theta(x)} = \frac{\theta(t)}{\theta(q(t))} = \frac{\theta(t)}{\theta(t+t^2)}.\eqno(8)$$
	As an example, let's take the right-hand limit of (8) as $t$ tends to 0. Since $t \rightarrow 0^+$ as $x \rightarrow 0^+$, we have
	$$\lim_{x \rightarrow 0^+} \frac{|BC|}{|ED|} = \lim_{t \rightarrow 0^+} \frac{\theta(t)}{\theta(t + t^2)} = \lim_{t \rightarrow 0^+} e^{-1/(t+1)} = e^{-1} \neq 1,\eqno(9)$$
	which finishes our construction.

	\section{Conclusions}~\

	As shown in Section 4, just using the intuitive perspective of geometry, we will not be able to explain why Arnold's proof fails at this point, since adding a $\theta(x)$ near 0 ti the graph makes no difference. Therefore, this geometric proof, altho not said to be wrong, does at least have its imperfections.\par
	If we look into the nature of the problem, we will find that the problem lies in the process of ``$BCED$ tends to a parallelogram''. If $f,\,g$ are analytic functions, then they approach each other at a rate close to the rate at which they approach $y = x$, and $BCED$ is indeed a parallelogram in the limit; whereas if $f,\,g$ are constructed as we have done in Section 4, $BCED$ only becomes more and more elongated and tends to ``diverge''. But the subtle differences in the properties of the functions, and the tendency of the graphs to change, cannot be captured in a static figure. Therefore, the discussion in this paper makes us realize once again that, altho intuitive graphical representations can help us quickly identify the essence of the problem, rigorous proofs are still indispensable. \par

\end{document}